\documentclass[12pt]{article}

\usepackage{amssymb,amsmath,exscale,times}

\usepackage{graphicx}
\usepackage{dsfont}
\usepackage[applemac]{inputenc}
\usepackage[T1]{fontenc}

\usepackage{color}
\definecolor{marin}{rgb}   {0.,   0.3,   0.7} 
\definecolor{rouge}{rgb}   {0.8,   0.,   0.} 
\definecolor{sepia}{rgb}   {0.8,   0.5,   0.} 
\usepackage[colorlinks,citecolor=marin,linkcolor=rouge,
            bookmarksopen,
            bookmarksnumbered
           ]{hyperref}

\newtheorem{lemma}{Lemma}[section]

\newtheorem{theorem}[lemma]{Theorem}
\newtheorem{proposition}[lemma]{Proposition}

\newtheorem{remark}[lemma]{Remark}
\newtheorem{example}[lemma]{Example}
\newtheorem{hypothesis}[lemma]{Hypothesis}
\newtheorem{notation}[lemma]{Notation}
\newtheorem{definition}[lemma]{Definition}
\newtheorem{conclusion}[lemma]{Conclusion}
\numberwithin{equation}{section}
\newcommand{\QED}{\mbox{}\hfill \raisebox{-0.2pt}{\rule{5.6pt}{6pt}\rule{0pt}{0pt}} 
          \medskip\par}             
\newenvironment{Proof}{\noindent
    \parindent=0pt\abovedisplayskip = 0.5\abovedisplayskip
    \belowdisplayskip=\abovedisplayskip{\bfseries Proof. }}{\QED}
\newenvironment{Proofof}[1]{\noindent
    \parindent=0pt\abovedisplayskip = 0.5\abovedisplayskip
    \belowdisplayskip=\abovedisplayskip{\bfseries Proof of #1. }}{\QED}

\newcommand{\Rc}{\mathsf{R}}

\newcommand{\Ac}{\mathcal{A}}
\newcommand{\ab}{\boldsymbol{a}}
\newcommand{\bb}{\boldsymbol{b}}

\newcommand{\dd}{\mathrm{d}}

\newcommand{\Hc}{\mathcal{H}}

\newcommand{\Ic}{\mathcal{I}}
\newcommand{\Jc}{\mathcal{J}}

\newcommand{\jb}{{\boldsymbol{j}}}
\newcommand{\ib}{{\boldsymbol{i}}}
\newcommand{\kb}{{\boldsymbol{k}}}
\newcommand{\Lc}{\mathcal{L}}

\newcommand{\ellb}{{\boldsymbol{\ell}}}
\newcommand{\N}{\mathbb{N}}
\newcommand{\Nc}{\mathcal{N}}

\newcommand{\Mc}{\mathcal{M}}
\newcommand{\Pc}{\mathcal{P}}

\newcommand{\R}{\mathbb{R}}

\newcommand{\C}{\mathbb{C}}

\newcommand{\T}{\mathbb{T}}
\newcommand{\Z}{\mathbb{Z}}

\newcommand{\Vc}{\mathcal{V}}

\newcommand{\W}{\mathcal{W}}

\newcommand{\Zc}{\mathcal{Z}}

\newcommand{\eps}{\varepsilon}

\newcommand{\va}[1]{|#1|}

\newcommand{\Norm}[2]{\|#1\|\left.\vphantom{T_{j_0}^0}\!\!\right._{#2}}         
             
\title{A Nekhoroshev type theorem for the nonlinear Schr\"odinger equation on the d-dimensional torus.}        
\author{Erwan Faou{\small$\,^1$} and Beno\^it Gr\'ebert{\small$\,^2$}\\[4ex]
{\small$\,^1$}\, \small INRIA \& ENS Cachan Bretagne,  \\[-1ex]
\small Avenue Robert Schumann F-35170 Bruz, France. \\[-1ex]
\it \small email: \tt Erwan.Faou@inria.fr\\[2ex]
{\small$\,^2$}\, \small Laboratoire de Math\'ematiques Jean Leray,
Universit\'e de Nantes,\\[-1ex]
\small 2, rue de la Houssini\`ere
F-44322 Nantes cedex 3, France. \\[-1ex]
\small\it email: \tt benoit.grebert@univ-nantes.fr
}       
\begin{document}
\maketitle
\abstract{We prove a Nekhoroshev type theorem for the nonlinear Schr\"odinger equation
$$
iu_t=-\Delta u+V\star u+\partial_{\bar u}g(u,\bar u)\ , \quad x\in \T^d,
$$
where $V$ is a typical smooth potential and $g$ is analytic in both variables.  More precisely we prove that if the initial datum is analytic in a strip of width $\rho>0$ with a bound on this strip equals to $\eps$ then, if $\eps$ is small enough, the solution of the nonlinear Schr\"odinger equation above remains analytic  in a strip of width $\rho/2$ and bounded on this strip by $C\eps$ during very long time of order $ \eps^{-\alpha|\ln \eps|^\beta}$ for some  constants $C> 0$, $\alpha>0$ and $\beta<1$.
\\[2ex]
}\\[2ex]
\tableofcontents
\section{Introduction and statements}
We consider the nonlinear Schr\"odinger equation
\begin{equation}\label{nls}
iu_t=-\Delta u+V\star u+\partial_{\bar u}g(u,\bar u)\ , \quad x\in \T^d, \quad t \in \R,
\end{equation}
where $V$ is a smooth convolution potential and $g$ is an analytic  function on a neighborhood of the origin in $\C^2$ which has a zero of order at least 3 at the origin and satisfies $g(z,\bar z)\in \R$. In more standard models, the convolution term is replaced by a multiplicative potential. The use of a convolution potential makes easier the analysis of the resonances.\\
For instance when  $g(u,\bar u)=\frac{a}{p+1}|u|^{2p+2}$ with $a\in \R$ and $p\in \N$,  we recover the standard NLS equation  
$iu_t=-\Delta u+V\star u+a|u|^{2p}u$. We notice that \eqref{nls} is a Hamiltonian system associated with the Hamiltonian function 
$$
H(u,\bar u ) = \int_{\T^d} \left(\va{\nabla u }^{2} +(V\star u
)\bar u  + 
g(u ,\bar u ) \right) \dd x.
$$
and the symplectic structure inherent to the complex structure, $i\dd u\wedge \dd\bar u$.\\

This equation has been considered with Hamiltonian tools  in two recent works.  In the first one (see \cite{BG03} and also \cite{BG06} and \cite{Bourgain} for related results) Bambusi \& Gr\'ebert prove a Birkhoff normal form theorem  adapted to this equation and obtain dynamical consequences on the long time behavior of the solutions with small initial Cauchy data in Sobolev spaces. More precisely they prove that if  the Sobolev norm of index $s$ of the initial datum $u_0$ is sufficiently small (of order $\varepsilon$) then the Sobolev norm of the solution is bounded by $2\varepsilon$ during very long time (of order $\varepsilon^{-r}$ with $r$ arbitrary). In the second one (see \cite{EK}) Eliasson \& Kuksin obtain a KAM theorem adapted to this equation. In particular they prove that, in a neighborhood of $u=0$, many of the invariant finite dimensional tori of the linear part of the equation  are preserved by small Hamiltonian perturbations. In other words, \eqref{nls} has many quasi-periodic solutions. In both cases non resonances conditions (not exactly the same) have to be imposed on the frequencies of the linear part and thus on the potential $V$. 

Both results are related to the stability of the zero solution which is an elliptic equilibrium of the linear equation. The first establishes the stability for polynomials times with respect to the size of the (small) initial  datum  while the second proves the stability for all time of certain solutions. In the present work we extend the technic of normal form and  we establish the stability for times of order $ \eps^{-\alpha|\ln \eps|^\beta}$ for some  constants $\alpha>0$ and $\beta<1$, $\eps$ being the size of the initial datum in an analytic space. \\

We now state precisely our result.
We assume that $V$ belongs to the following space ($m> d/2$, $R>0$)
\begin{equation}\label{pot}
\W_{m}=\{ V(x)=\sum_{a\in \Z^d}v_a
e^{i a\cdot x}
\mid v'_{a}:={v_a}{(1+\va{a})^m}/R \in [-1/2,1/2]
\mbox{ for any }a\in \Z^d \}
\end{equation}
that we endow with the product probability measure. Here, for $a = (a_1,\ldots,a_d) \in \Z^d$, $|a|^2 = a_1^2 + \cdots + a_d^2$. 

For $\rho>0$, we denote by $\Ac_\rho\equiv \Ac_\rho(\T^d ; \C)$ the
space of functions $\phi$ that are analytic on the complex neighborhood of $d$-dimensional torus $\T^d$ given by $I_\rho=\{x+iy \mid x\in \T^d,\ y\in \R^d \mbox{ and }\ |y|< \rho\}$ and continuous on the closure of this strip. We then denote by
 $|\cdot|_\rho$ the usual norm on $\Ac_\rho$
$$
|\phi|_\rho=\sup_{z\in I_\rho}|\phi(z)|.
$$ 
We note that $(\Ac_\rho,|\cdot|_\rho)$ is a Banach space.\\
Our main result is a Nekhoroshev type theorem:
\begin{theorem}\label{main} 
There exists a subset   $\Vc \subset \W_{m}$  of full measure, such that  for $V \in \Vc$,  $\beta<1$ and $\rho >0$, the following holds:  there exist $C>0$ and  $\varepsilon_0>0$  such that
if 
$$u_0\in \Ac_{2\rho} \quad \mbox{ and  }\quad |u_0|_{2\rho} =\eps\leq \eps_0$$
 then the solution of \eqref{nls} with initial datum $u_0$ exists for times $|t|\leq \eps^{-\alpha(\ln \eps)^\beta}$ and satisfies
\begin{equation}\label{mainest1}
|u(t)|_{\rho/2}\leq C \eps \quad \mbox{ for }\quad |t|\leq  \eps^{-\sigma_\rho |\ln \eps|^\beta},
\end{equation}
with  $\sigma_\rho=\min\{\frac 1 8,\frac \rho 2\}$.\\
Furthermore, writing $u(t)=\sum_{k\in \Z^d}\xi_k(t)e^{ik\cdot x}$, we have 
\begin{equation}\label{mainest2}
\sum_{k\in\Z^d}e^{ \rho |k|  } \big | |\xi_k(t)|-|\xi_k(0)|\big|\leq \eps^{3/2} \quad \mbox{ for }\quad  |t|\leq  \eps^{-\sigma_\rho|\ln \eps|^\beta}.
\end{equation} 
\end{theorem}
Estimate \eqref{mainest2} asserts that there is almost no variation of the actions\footnote{Here the actions are the modulus of the Fourier coefficients to the square, $I_k=|\xi_k|^2$.} and in particular no possibility of weak turbulence, i.e. exchanges between low Fourier modes and high Fourier modes. This kind of turbulence may induce the growth of the Sobolev  norm $\sum (1+|k|^s)^2|\xi_k|^2$ ($s>1$) of the solution as recently proved in \cite{Iteam}.\\
In finite dimension $n$, the standard Nekhoroshev result \cite{Nek77} controls the dynamic over times of order $\exp \left(\frac{-\alpha}{\eps^{1/(\tau +1)}}\right)$ for some $\alpha>0$ and $\tau>n+1$   (see for instance \cite{ben85, GG85, P93}) which is of course much better than  $ \eps^{-\alpha|\ln \eps|^\beta}=e^{-\alpha|\ln \eps|^{(1+\beta)}}$. Nevertheless this standard result does not extend to the infinite dimensional context. Actually, when $n\to \infty$, that $\eps^{-1/(\tau +1)}$ can be transformed in $|\ln \eps|^{(1+\beta)}$ is a good news!\\
The only previous work in the direction of Nekhoroshev estimates for PDEs was obtained by Bambusi in \cite{Bam99}. He also worked in spaces of analytic functions in a strip and for times of order $e^{-\alpha |\ln \eps|^{1+\beta}}$, nevertheless the control of the solution was not obtained uniformly in a strip but  in a complicated way involving the Fourier coefficients of the solution.\\
 
We now focus on the three main differences with the previous works on normal forms:
\begin{itemize}
\item we crucially use the zero momentum condition: in the Fourier space, the nonlinear term contains only monomials $z_{j_1}\cdots z_{j_k}$ with $j_1+\cdots  +j_k=0$ (cf. Definition \ref{zero}). This property allows to control the largest index by the others.

\item we use  $\ell^1$-type norms   to control the Fourier coefficients and the vector fields instead of $\ell^2$-type norms as usual. Of course this choice does not allow to work in Hilbert spaces and makes obligatory a slight lost of regularity each time the estimates are transposed from the Fourier space to the initial space of analytic functions. But it turns out that this choice makes much more simpler the estimates on the vector fields (cf. Proposition \ref{P1}  below and \cite{FG10} for a similar framework in the context of numerical analysis). 

\item we notice that the vector field of a monomial,$z_{j_1}\cdots z_{j_k}$ containing at least three Fourier modes $z_\ell$ with large indices $\ell$ induces a flow whose dynamics is under control during very long time in the sense that the dynamic almost excludes exchanges between high Fourier modes and low Fourier modes  (see Proposition \ref{Pcrux}). In \cite{Bam03} or \cite{BG06}, such terms were neglected since the vector field of a monomial containing at least three Fourier modes with large indices is small in \textit{Sobolev norm} (but not in analytic norm) and thus will almost keep invariant all the modes. This more subtle analysis for monomials  was still used in \cite{FGP2}.

\end{itemize}
Finally we notice that our method could be generalized by considering not only zero momentum monomials but also monomials with   finite or exponentially decreasing momentum. This would   certainly allow to consider a nonlinear Schr\"odinger equation with a multiplicative potential $V$ and nonlinearities depending periodically on $x$: 
$$
iu_t=-\Delta u+V u+\partial_{\bar u}g(x,u,\bar u)\ , \quad x\in \T^d.
$$
Nevertheless this generalization would generate a lot of technicalities and we prefer to focus here on the simplicity of the arguments.

\section{Setting and Hypothesis}
\subsection{Hamiltonian formalism}
The equation  \eqref{nls} is a semi linear PDE locally well posed in the Sobolev space $H^2(\T^d)$ (see for instance \cite{Cazenave}).
Let $u$ be a (local) solution of \eqref{nls} and consider $(\xi,\eta) = (\xi_a,\eta_a)_{a \in \Z^d}$ the Fourier coefficients of $u$, $\bar u$ respectively, i.e. 
\begin{equation}
\label{Exieta}
u (x)= \sum_{a\in \Z^d} \xi_{a} 
e^{i a\cdot x} \quad\mbox{and}\quad
 \bar u (x)= \sum_{a\in \Z^d}
\eta_{a} 
e^{-i a\cdot x}.
\end{equation}
A standard calculus shows that $u$ is solution in $H^2(\T^d)$ of \eqref{nls} if and only if $(\xi,\eta)$ is a solution in\footnote{As usual, $\ell^2_2=\{(\xi_a)_{a\in\Z^d}\mid \sum (1+|a|^2)|\xi_a|^2<+\infty\}.$} $\ell^2_2\times \ell^2_2$ of the system
\begin{equation}\label{nlsxieta}
\left\{\begin{array}{rll} \displaystyle \dot \xi_a&= -i\omega_a \xi_a -i\frac{\partial P}{\partial \eta_a},  \quad &a\in \Z^d,
 \\[2ex]
  \dot \eta_a&= i\omega_a \eta_a -i\frac{\partial P}{\partial \xi_a},  \quad &a\in \Z^d,
 \end{array}
 \right.
\end{equation}
where the linear frequencies are given by
$\omega_{a}= |a|^{2} +  v_a$ where as in \eqref{pot}, $V=\sum v_a e^{ia\cdot x}$, and the nonlinear part is given by
\begin{equation}\label{P}
P(\xi,\eta)=  \frac{1}{(2\pi)^d}  \int_{\T^d} g(\sum \xi_{a} 
e^{i a\cdot x} , \sum
\eta_{a} 
e^{-i a\cdot x} ) \,  \dd x.
\end{equation}
This system is  reinterpreted in a Hamiltonian context endowing the set of couples $(\xi_a,\eta_a) \in \C^{\Z^d} \times \C^{\Z^d}$  with the symplectic structure
\begin{equation}
\label{Esymp}
i \sum_{a \in \Z^d} \dd \xi_a \wedge \dd \eta_a. 
\end{equation}
We define the set $\Zc = \Z^d \times \{ \pm 1\}$. 
For $j = (a,\delta) \in \Zc$, we define $|j| = |a|$ and we denote by $\overline{j}$ the index $(a,-\delta)$. \\
We identify a couple $(\xi,\eta)\in \C^{\Z^d} \times \C^{\Z^d}$ with 
$(z_j)_{j \in \Zc} \in \C^{\Zc}$ via the formula
\begin{equation}
\label{Ezj}
j = (a,\delta) \in \Zc  \Longrightarrow 
\left\{
\begin{array}{rcll}
z_{j} &=& \xi_{a}& \mbox{if}\quad \delta = 1,\\[1ex]
z_j &=& \eta_a & \mbox{if}\quad \delta = - 1.
\end{array}
\right.
\end{equation}
By  a slight abuse of notation, we often write $z = (\xi,\eta)$ to denote such an element. \\
For a given $\rho > 0$, we consider the Banach space 
$\Lc_\rho$ made of elements $z \in \C^\Zc$ such that 
$$
\Norm{z}{\rho} := \sum_{j \in \Zc} e^{\rho|j|} |z_j| < \infty,
$$
and equipped with the symplectic form \eqref{Esymp}. We say that $z\in \Lc_\rho$ is {\em real} when $z_{\overline{j}} = \overline{z_j}$ for any $j\in \Zc$. In this case, we write $z=(\xi,\bar\xi)$ for some $\xi\in \C^{\Z^d}$. In this situation, we can associate with $z$ the function $u$ defined by \eqref{Exieta}. \\ 
The next lemma shows the relation with the space $\Ac_\rho$ defined above: 
\begin{lemma}\label{AL}
Let $u$ be a complex valued function analytic on a neighborhood of  $\T^d$, and let $(z_j)_{j \in \Zc}$ be the sequence of its Fourier coefficients defined by \eqref{Exieta} and \eqref{Ezj}. Then for all $\mu < \rho$, we have 
\begin{align}
\label{Ediff}
\text{if}\quad  u\in \Ac_\rho\quad & \text{then}\quad z\in \Lc_\mu  \quad \text{and} \quad \Norm{z}{\mu} \leq c_{\rho,\mu} |u|_\rho \ ;\\ \label{Ediff2}
\text{if}\quad z\in \Lc_\rho\quad    &\text{then}\quad u\in \Ac_\mu \quad  \text{and} \quad \va{u}_{\mu} \leq c_{\rho,\mu} \Norm{z}{\rho},
\end{align}
where $c_{\rho,\mu}$ is a constant depending on $\rho$ and $\mu$ and the dimension $d$. 
\end{lemma}
\begin{Proof}
Assume that $u \in \Ac_\rho$. Then by Cauchy formula, we have for all $j\in\Zc$,  $|z_j| \leq |u|_\rho e^{-\rho |j| }$. Hence for $\mu < \rho$, we have 
$$
\Norm{z}{\mu} \leq |\phi|_\rho\sum_{j \in \Zc} e^{(\mu-\rho) |j| } \leq  |\phi|_\rho \Big(2 \sum_{n \in \Z} e^{\frac{(\mu-\rho)}{\sqrt{d}}|n| }  \Big)^d   \leq \left( \frac{2}{1 - e^{\frac{(\mu-\rho)}{\sqrt{d}}}} \right)^d|u|_\rho.  
$$
Conversely, assume that $z \in \Lc_\rho$. Then  $|\xi_a| \leq \Norm{z}{\mu} e^{-\rho |a|}$ for all $a \in \Z^d$, and thus by \eqref{Exieta},  we get for all $x\in \T^d$ and $y\in \R^d$ with $|y|\leq \mu$,
$$|u(x+iy)|\leq\sum_{  a \in \Z^d   } |\xi_a|e^{|ay|}\leq  \Norm{z}{\rho}\sum_{  a \in \Z^d   } e^{-(\rho-\mu) |a|}\leq    \left( \frac{2}{1 - e^{\frac{(\mu-\rho)}{\sqrt{d}}}} \right)^d    \Norm{z}{\rho}.$$
Hence $u$  is bounded on the strip $I_\mu$. 
\end{Proof}

For a function $F$ of $\mathcal{C}^1(\Lc_\rho,\C)$, we define its Hamiltonian vector field by $X_F=J\nabla F$ where  $J$ is the symplectic operator on $\Lc_\rho$ induced  by the symplectic form \eqref{Esymp},
$
\nabla F(z) = \left( \frac{\partial F}{\partial z_j}\right)_{j \in \Zc}$
and where by definition we set for $j = (a,\delta) \in \Z^d \times \{ \pm 1\}$, 
$$
 \frac{\partial F}{\partial z_j} =
  \left\{\begin{array}{rll}
 \displaystyle  \frac{\partial F}{\partial \xi_a} & \mbox{if}\quad\delta = 1,\\[2ex]
 \displaystyle \frac{\partial F}{\partial \eta_a} & \mbox{if}\quad\delta = - 1.
 \end{array}
 \right.
$$
For two functions $F$ and $G$, the Poisson Bracket is (formally) defined as
\begin{equation}\label{poisson}
\{F,G\} = \nabla F^T J \nabla G = i \sum_{a \in \Z^d} \frac{\partial F}{\partial \eta_a}\frac{\partial G}{\partial \xi_a} -  \frac{\partial F}{\partial \xi_a}\frac{\partial G}{\partial \eta_a}.  
\end{equation}
We say that a Hamiltonian function $H$ is 
 {\em real } if $H(z)$ is real for all real $z$. 
 
\begin{definition}\label{def:2.1}
 For a given $\rho > 0$, we  denote by $\Hc_\rho$ the space of real Hamiltonians $P$ satisfying 
$$
P \in \mathcal{C}^{1}(\Lc_\rho,\C), \quad \mbox{and}\quad 
X_P \in \mathcal{C}^{1}(\Lc_\rho,\Lc_\rho). 
$$
\end{definition}
Notice that for $F$ and $G$ in $\Hc_\rho$ the formula \eqref{poisson} is well defined. 
With a given Hamiltonian function $H \in \Hc_\rho$, we associate the Hamiltonian system
$$
\dot z =   X_H(z) =    J \nabla H(z)
$$
which also reads
\begin{equation}\label{Eham2}
\displaystyle \dot \xi_a - -i\frac{\partial H}{\partial \eta_a}  \quad \mbox{and} \quad 
  \dot \eta_a= i\frac{\partial H}{\partial \xi_a},  \quad a\in \Z^d.
\end{equation}
We define the local flow $\Phi_H^t(z)$ associated with the previous system (for an interval of times $t \geq 0$ depending a priori on the initial condition $z$). Note that if $z = (\xi,\bar \xi)$ and if  $H$ is real, the flow $(\xi^t,\eta^t) = \Phi_H^t(z)$ is also real, $\xi^t = \bar {\eta}^t$ for all $t$. Further, choosing the Hamiltonian given by
$$
H(\xi,\eta)=\sum_{a\in \Z^d}\omega_a \xi_a\eta_a + P(\xi,\eta),$$
$P$ being given by \eqref{P},
we recover the system \eqref{nlsxieta}, i.e.  the expression of the NLS equation \eqref{nls} in Fourier modes. 
\begin{remark}\label{remH0}
The quadratic Hamiltonian $H_0=\sum_{a\in \Z^d}\omega_a \xi_a\eta_a$ corresponding to the linear part of \eqref{nls} does not belong to $\Hc_\rho$. Nevertheless it generates a flow which maps $\Lc_\rho$ into $\Lc_\rho$ explicitly given for all time $t$ and for all indices $a$ by  $\xi_a(t)=e^{-i\omega_a t}\xi_k(0)$, $\eta_a(t)=e^{i\omega_a t}\eta_k(0)$. On the contrary, we will see that, in our setting, the nonlinearity $P$
belongs to $\Hc_\rho$. \end{remark}

\subsection{Space of polynomials}
In this subsection we  define a class of polynomials  on $\C^\Zc$.\\
We first need more notations concerning multi-indices: let $\ell\geq 2$ and $\jb = (j_1,\ldots,j_\ell)\in \Zc^\ell$ with $j_i=(a_i,\delta_i)$, we define
\begin{itemize}
\item the monomial associated with $\jb$ :
$$\quad z_\jb = z_{j_1}\cdots z_{j_\ell},$$
\item the momentum of $\jb$ :
\begin{equation}
\label{EMcb}
\Mc(\jb) = a_{1} \delta_{1} + \cdots + a_\ell \delta_\ell,
\end{equation}
\item the divisor associated with $\jb$ :
\begin{equation}
\label{EOmega}
\Omega(\jb) = 
\delta_1\omega_{a_1} + \cdots  + \delta_\ell\omega_{a_\ell}
\end{equation}
where, for $a\in \Z^d$,  $\omega_a=|a|^2+v_a$ are the frequencies of the linear part of \eqref{nls}. 
\end{itemize}
 We then define the set of indices with {\bf zero momentum} by
\begin{equation}
\label{EIr}
 \Ic_\ell =  \{  \jb = (j_1,\ldots,j_\ell) \in \Zc^{\ell}, \quad \mbox{with}\quad \Mc(\jb) = 0\}.
\end{equation}
On the other hand, we say that $\jb = (j_1,\ldots,j_\ell) \in\Zc^\ell$ is {\bf resonant}, and we write $\jb \in \Nc_r$, if $\ell$ is even and $\jb=\ib \cup \bar \ib$ for some choice of   $\ib\in \Zc^{\ell/2}$.   
In particular,  if  $\jb$ is resonant then its associated divisor vanishes, $\Omega(\jb)=0$, and its associated monomials depends only on the actions: 
$$
z_\jb = z_{j_1}\cdots z_{j_r} = \xi_{a_1}\eta_{a_1} \cdots \xi_{a_{\ell/2}} \eta_{a_{\ell/2}}
= I_{a_1} \cdots I_{a_{\ell/2}}, 
$$
where for all $a \in \Z^d$, 
$
I_{a}(z) = \xi_a \eta_a
$
denotes the action associated with the index $a$.\\ 
Finally we note that if $z$ is real, then  $I_a(z) = |\xi_a|^2$ and we remark that for odd $r$ the resonant set  $\Nc_r$ is the empty set.  

\begin{definition}\label{zero}
Let $k\geq 2$, a (formal) polynomial $P(z) = \sum a_{\jb} z_{\jb} $ belongs to $ \Pc_k$ if $P$ is real, of  degree $k$, have a zero of order at least $2$ in $z = 0$, and if 
\begin{itemize}
\item $P$ contains only monomials  having zero momentum, i.e.\ such that $\Mc(\jb)=0$ when $a_\jb \neq 0$ and thus $P$  reads
\begin{equation}
\label{EexpP}
P(z) = \sum_{\ell = 2}^k \sum_{\jb \in \Ic_\ell} a_{\jb} z_{\jb}
\end{equation}
with the relation $a_{\bar\jb} = \bar{a}_{\jb}$. 
\item The coefficients $a_{\jb}$ are bounded, i.e.\ $\forall\, \ell = 2,\ldots,k$, $\displaystyle \sup_{\jb \in \Ic_\ell} |a_\jb| < + \infty. 
$
\end{itemize}
\end{definition}
We endow $ \Pc_k$ with the norm
\begin{equation}
\label{EnormP}
\Norm{P}{} = \sum_{\ell = 2}^k \sup_{\jb \in \Ic_\ell} |a_\jb|. 
\end{equation} 
The zero momentum assumption  in Definition \ref{zero} is crucial  to obtain the following Proposition:
\begin{proposition}
\label{P1}
Let $k \geq 2$  and $\rho > 0$. We have  $\Pc_k \subset \Hc_\rho$,  and for $P$ a homogeneous polynomial of degree $k$ in  $ \Pc_k$,  we have the estimates
\begin{equation}
\label{Epot}
|P(z)| \leq \Norm{P}{}\Norm{z}{\rho}^k 
\end{equation}
and 
\begin{equation}
\label{Echamp}
\forall\, z \in \Lc_\rho, \quad   \Norm{X_P(z) }{\rho} \leq 2 k \Norm{P}{}\Norm{z}{\rho}^{k-1}. 
\end{equation}
Eventually,  for $P\in \Pc_k$ and $Q \in \Pc_\ell$, then $\{P,Q\} \in \Pc_{k+ \ell - 2}$ and we have the estimate
\begin{equation}
\label{Ebrack}
\Norm{\{P,Q\}}{} \leq 2 k\ell\Norm{P}{}\Norm{Q}{}.
\end{equation}
\end{proposition}

\begin{Proof}
Let
$$
P(z) = \sum_{\jb \in \Ic_k} a_\jb z_\jb, 
$$
we have
$$
|P(z)|\leq   \Norm{P}{}\sum_{\jb \in \Zc^k}  |z_{j_1}|\cdots |z_{j_k}|\leq \Norm{P}{} \Norm{z}{\ell^1}^k  \leq  \Norm{P}{} \Norm{z}{\rho}^k 
$$
and the first inequality \eqref{Epot} is proved. \\
To prove the second estimate,  let us take $\ell \in \Zc$ and  calculate using the zero momentum condition,
$$
\left|\frac{\partial P}{\partial z_\ell} \right| \leq k\Norm{P}{}\sum_{ \substack{\jb \in \Zc^{k-1} \\ \Mc(\jb) = -  \Mc(\ell) }}|z_{j_1} \cdots z_{j_{k-1}}|.
$$
 Therefore
$$
\Norm{X_P(z) }{\rho}  =\sum_{\ell\in\Zc}e^{\rho|\ell|} \left|\frac{\partial P}{\partial z_\ell} \right|  \leq k \Norm{P}{}\sum_{\ell\in\Zc} \sum_{ \substack{\jb \in \Zc^{k-1} \\ \Mc(\jb) = -  \Mc(\ell) }}  e^{\rho|\ell|} |z_{j_1} \cdots z_{j_{k-1}}|.
$$
But if $\Mc(\jb) = - \Mc(\ell)$, 
$$
e^{\rho|\ell|} \leq \exp\big(\rho( |j_1| + \cdots + |j_{k-1}|)\big) \leq \prod_{n = 1,\ldots,k-1} e^{\rho |j_n|}. 
$$
Hence, after summing in $\ell$ we get\footnote{Take care that $\Mc(a,\delta) =  \Mc(-a,-\delta)$ whence the coefficient 2.}
$$
\Norm{X_{P}(z) }{\rho}\leq 2 k \Norm{P}{} \sum_{\jb \in \Zc^{k-1}} e^{\rho |j_1|}|z_{j_1}| \cdots  e^{\rho |j_{k-1}|}| z_{j_{k-1}}| 
\leq 2 k  \Norm{P}{}\Norm{z}{\rho}^{k-1}
$$
which yields \eqref{Echamp}. 

Assume now that $P$ and $Q$ are homogeneous polynomials of degrees $k$ and $\ell$ respectively and with coefficients $a_\kb$, $\kb \in \Ic_k$ and $b_\ellb$, $\ellb \in \Ic_\ell$. It is clear that $\{P,Q\}$ is a monomial of degree $k + \ell - 2$ satisfying the zero momentum condition. Furthermore writing
$$
\{P,Q\}(z) = \sum_{ \jb \in \Ic_{k+\ell-2} } c_{\jb} z_\jb,
$$
$c_\jb$ expresses  as a sum of coefficients $a_{\kb}b_\ellb$ for which there exists an $a \in \Z^d$ and $\epsilon \in \{\pm 1\}$ such that 
$$
(a,\epsilon) \subset \kb \in \Ic_k \quad \mbox{and} \quad (a,-\epsilon) \subset \ellb \in \Ic_\ell, 
$$
and such that if for instance $(a,\epsilon) = k_1$ and $(a,-\epsilon) = \ell_1$, we necessarily have  $(k_2,\ldots,k_k,\ell_2,\ldots,\ell_\ell) = \jb$. Hence for a given $\jb$, the zero momentum  condition on $\kb$ and on $\ellb$ determines  the value of $\epsilon a$ which in turn determines two possible value of $(\epsilon,a)$.\\
This proves \eqref{Ebrack} for monomials.  The extension to polynomials follows from the definition of the norm \eqref{EnormP}. \\ 
The last assertion, as well as the fact that the Poisson bracket of two real Hamiltonian is real, immediately  follow  from the definitions.
\end{Proof}

\subsection{Nonlinearity}\label{sec:NL}

The nonlinearity $g$ in \eqref{nls} is assumed to be complex analytic in a neighborhood of $\{0,0\}$ in $\C^2$. So there exist positive constants $M$ and $R_0$ such that  the Taylor expansion 
$$
g(v_1,v_2) = \sum_{k_1,k_2\geq 0} \frac{1}{k_1! k_2!}\partial_{k_1} \partial_{k_2} g (0,0) v_1^k v_2^\ell
$$
is uniformly convergent and bounded by $M$ on the ball   $|v_1|+|v_2|\leq2R_0$ of $\C^2$. Hence, formula \eqref{P} defines an analytic function on the ball $\Norm{z}{\rho}\leq R_0$ of $\Lc_\rho$ and we have
$$
P(z) =
 \sum_{k \geq 0} P_k(z) 
$$
where, for all $k\geq 0$, $P_k$ is a homogeneous polynomial defined by
$$
P_k =  \sum_{k_1 + k_2 = k} 
\sum_{(\ab,\bb) \in (\Z^d)^{k_1}  \times (\Z^d)^{k_1}} p_{\ab,\bb} \xi_{a_1}\cdots \xi_{a_{k_1}} \eta_{b_1}\cdots \eta_{b_{k_2}}
$$
with
$$
p_{\ab,\bb} = 
 \frac{1}{k_1! k_2!}\partial_{k_1} \partial_{k_2} g (0,0) \int_{\T^d} e^{i \Mc(\ab,\bb)\cdot x} \, \dd x,
$$
and $\Mc(\ab,\bb) = a_{1} + \cdots +a_{k_1} - b_{1} - \cdots- b_{k_2}$ is the moment of $\xi_{a_1}\cdots \xi_{a_{k_1}} \eta_{b_1}\cdots \eta_{b_{k_2}}$. Therefore it is clear that $P_k$ satisfies the zero momentum condition and thus
 $P_k \in \Pc_k$ for all $k\geq 0$. Furthermore  $\Norm{P_k}{} \leq M R_0^{-k}$. 

\subsection{Non resonance condition}\label{sec:NR}
In order to control the divisors \eqref{EOmega}, we need to impose a non resonance condition on the linear frequencies $\omega_a$, $a\in \Z^d$.\\ 
For $r \geq 3$ and $\jb = (j_1,\ldots,j_r) \in \Zc^r$, we define $\mu(\jb)$ as the third largest integer amongst $|j_1|, \cdots, |j_r|$ and we recall that $\jb\in \Zc^r$ is said resonant if $r$ is even and $\jb=\ib\cup\bar\ib$ for some $\ib\in\Zc^{r/2}$.

\begin{hypothesis}
\label{NR} 
There exist $\gamma > 0$,   $\nu > 0$   and   $c_0 > 0$    such that for all $r \geq 3$ and all $\jb \in \Zc^r$ non resonant, we have
\begin{equation}
\label{Enonres1}
|\Omega(\jb)| \geq  \frac{  \gamma c_0^r   }{ \mu(\jb)^{\nu r}}.
\end{equation}
\end{hypothesis}
Recall that for $V = \sum_{a \in \Z^d} v_a e^{ia\cdot x}$ in the  space $\W_{m}$ defined in \eqref{pot}, the frequencies read
$$
\omega_{a}=|a|^2 + v_a = |a|^2+\frac{Rv_{a}'}{(1+|a|)^m},\quad a \in \Z^d.
$$
In Appendix we prove
\begin{proposition} \label{res.1}
Fix $\gamma >0$ small enough and $m>d/2$.  There exist positive
constants    $c_0$    and $\nu$ depending only on $m$, $R$ and $d$,  and a set $F_{\gamma} \subset \W_{m}$
whose measure is larger than $1- 4 \gamma^{1/7}$ such that if $V\in F_{\gamma}$
then \eqref{Enonres1} holds true for all non resonant $\jb\in\Zc^r$ and all $r\geq 3$.
\end{proposition}
Thus Hypothesis \ref{NR} is satisfied for all $V\in \Vc$ where
\begin{equation}\label{V}
\Vc=\cup_{\gamma>0}F_{\gamma}
\end{equation}
is a subset of full measure in $\W_m$.
\subsection{Normal forms}

We fix an index $N \geq 1$. For a fixed integer $k \geq 3$, we set 
$$
\Jc_k(N) = \{ \jb \in \Ic_k\, | \, \mu(\jb) > N \}. 
$$

\begin{definition} Let $N$ be an integer.
We say that a polynomial $Z \in \Pc_k$ is in $N$-normal form if it can be written
$$
Z = \sum_{\ell = 3}^k\ \ \sum_{\jb \in \Nc_\ell \cup\Jc_\ell(N)   }a_\jb z_\jb
$$
In other words, $Z$ contains either monomials depending only of the actions or monomials whose indices $\jb$ satisfies $\mu(\jb) > N$, i.e. monomials involving at least three modes with index greater than N. 
\end{definition}
We now motivate the introduction of such normal form. First, we recall the 
\begin{lemma}
\label{Lcomp}
let $f:\R \to \R_+$ a continuous function, and $y:\R \to \R_+$ a differentiable function satisfying the inequality
$$
\forall\, t \in \R, \quad  \frac{\dd}{\dd t} y(t) \leq 2 f(t) \sqrt{y(t)}.
$$
Then we have the estimate
$$
\forall\, t \in \R, \quad \sqrt{y(t)} \leq \sqrt{y(0)} + \int_0^t f(s) \, \dd s.
$$
\end{lemma}
\begin{Proof} Let  $\epsilon>0$ and define $y_\epsilon = y+\epsilon$ which is a non negative function whose square root is derivable. We have
$$
 \frac{\dd}{\dd t}\sqrt{ y_\epsilon(t)} \leq 2 f(t) \frac{\sqrt{y(t)}}{\sqrt{ y_\epsilon(t)}}\leq 2 f(t)
$$
and thus 
$$
\sqrt{y_\eps(t)} \leq \sqrt{y_\epsilon(0)} + \int_0^t f(s) \, \dd s.
$$
The claim is obtained when $\epsilon\to 0$.
\end{Proof}
For a given number $N$ and $z \in \Lc_\rho$ we define 
$$ 
\Rc^N_{\rho}(z) = \sum_{|j| > N} e^{\rho|j|} |z_j|.
$$ 
Notice that if $z \in \Lc_{\rho+\mu}$ then
\begin{equation}
\label{RN}
\Rc^N_{\rho}(z) \leq e^{-\mu N}\Norm{z}{\rho+\mu}.
\end{equation}
\begin{proposition}
\label{Pcrux}
Let $N \in \N$ and $k \geq 3$. 
 Let $Z$ a homogeneous polynomial of degree $k$ in $N$-normal form. Let $z(t)$ be a real solution of the flow associated with the Hamiltonian $H_0+Z$. Then we have 
\begin{equation}
\label{ERN}
\Rc_\rho^N(t) \leq \Rc_\rho^N(0) +   4    k^3 \Norm{Z}{}\int_0^t\Rc_\rho^N(s)^2 \Norm{z(s)}{\rho}^{k-3}  \dd s
\end{equation}
and
\begin{equation}
\label{ENN}
\Norm{z(t)}{\rho} \leq \Norm{z(0)}{\rho}  +   4    k^3 \Norm{Z}{}\int_0^t\Rc_\rho^N(s)^2 \Norm{z(s)}{\rho}^{k-3}  \dd s
\end{equation}
\end{proposition}
\begin{Proof}
Let $a \in \Z^d$ be fixed, and let $I_a(t) = \xi_a(t) \eta_a(t)$ the actions associated with the solution of the Hamiltonian system induced by $H_0+Z$. 
We have using \eqref{Ebrack} and $H_0=H_0(I)$,
$$
| e^{2 \rho|a|} \dot I_a |= |e^{2 \rho|a|} \{I_a,Z\} |
\leq 2 k \Norm{Z}{}|e^{\rho|a|} \sqrt{I_a} | \Big(\sum_{\substack{\Mc(\jb) = \pm a\\   \text{2 indices}>  N  }} e^{\rho|a|}|z_{j_1} \cdots z_{j_{k-1}}|\Big)
$$
Using the previous Lemma, we get
\begin{equation}
\label{marre}
e^{\rho|a|} \sqrt{I_a(t)}  \leq  e^{\rho|a|} \sqrt{I_a(0)}  + 
 2 k \Norm{Z}{} \int_0^t \Big(\sum_{\substack{\Mc(\jb) = \pm a\\   \text{2 indices}> N  }} e^{\rho|j_1|}|z_{j_1}| \cdots e^{\rho|j_{k-1}|}|z_{j_{k-1}}|\Big) \dd s.
\end{equation}
Ordering the multi-indices in such way $|j_1|$ and $|j_2|$ are the largest,  and using the fact that $z(t)$ is real (and thus $|z_j| = \sqrt{I_a}$ for $j = (a, \pm 1) \in \Zc$), we obtain after summation in $|a| >  N$
\begin{align*}
  \Rc_\rho^N(z(t)) &\leq \Rc_\rho^N(z(0)) +     4  k^3 \Norm{Z}{} \int_0^t \Big(\sum_{\substack{|j_1|,|j_2|\geq N\\ j_3,\cdots, j_k\in \Zc}} e^{\rho|j_1|}|z_{j_1}| \cdots e^{\rho|j_{k-1}|}|z_{j_{k-1}}|\Big) \dd s\\
  &\leq \Rc_\rho^N(0) +   4   k^3 \Norm{Z}{}\int_0^t\Rc_\rho^N(s)^2 \Norm{z(s)}{\rho}^{k-3}  \dd s.
\end{align*}
In the same way we obtain \eqref{ENN}. 
\end{Proof}
\begin{remark}
\label{choixN}
These estimates will be crucially used in the final bootstrap argument. In particular, along the solution associated with a Hamiltonian in $N$-normal formal and initial datum 
$\Norm{z_0}{2\rho}=\varepsilon$. Then as $\Rc_\rho^N(z_0) = \mathcal{O}(\varepsilon e^{-\rho N})$,  Eqns. \eqref{ERN}-\eqref{ENN} guarantee that $\Rc_\rho^N(z(t))$ remains of order  $\mathcal{O}(\varepsilon e^{-\rho N})$ and the norm of $z(t)$ remains of order $\varepsilon$  over exponentially long time $ t = \mathcal{O}(e^{\rho N})$. 
\end{remark}
%
The next result is an easy consequence of the non resonance condition and the definition of the normal forms: 
\begin{proposition}
\label{Phomolo}
Assume that the non resonance condition \eqref{Enonres1} is satisfied, and let $N$ be fixed. Let $Q$ be a homegenous polynomial of degree $k$. Then the homological equation 
\begin{equation}
\label{Ehomo}
\{ \chi,H_0 \} - Z = Q
\end{equation}
admits a polynomial solution $(\chi,Z)$ homogeneous of degree $k$, such that $Z$ is in $N$-normal form, and such that 
\begin{equation}
\label{Ehomol}
\Norm{Z}{} \leq  \Norm{Q}{} \quad \mbox{and}\quad \Norm{\chi}{} \leq  \frac{N^{\nu k}}{\gamma c_0^k}\Norm{Q}{}
\end{equation}
\end{proposition}
\begin{Proof}
Assume that $Q = \sum_{\jb \in \Ic_k} Q_\jb z_\jb$ and search $Z = \sum_{\jb \in \Ic_k} Z_\jb z_\jb$ and $\chi = \sum_{\jb \in \Ic_k} \chi_\jb z_\jb$ such that \eqref{Ehomo} be satisfied. Then the equation \eqref{Ehomo} can be written in term of polynomial coefficients
$$
i \Omega(\jb) \chi_\jb - Z_\jb = Q_\jb, \quad j\in  \Ic_k,   
$$
where $\Omega(\jb)$ is defined in \eqref{EOmega}. We then define 
\begin{equation*}\begin{array}{llll}
&Z_\jb = Q_\jb\quad& \text{and}\quad \chi_{\jb} = 0 &\text{if} \quad \jb \notin \Nc_k \text{ or } \mu(\jb)\leq N,\\
&Z_\jb = 0\quad &\text{and}\quad \chi_{\jb} = \frac{Q_\jb}{i \Omega(\jb)}  &\text{if} \quad \jb \in \Nc_k \text{ and } \mu(\jb)> N.
\end{array}\end{equation*}
In view of \eqref{Enonres1}, this yields \eqref{Ehomol}. 
\end{Proof}

\section{Proof of  the main Theorem}
\subsection{Recursive equation}

We aim at constructing a canonical transformation $\tau$ such that in the new variables, the Hamiltonian $H_0 + P$ is under normal form modulo a small remainder term. 
Using Lie transforms to generate $\tau$, the problem can be written: 
Find polynomials  $\chi = \sum^r_{k = 3} \chi_{k}$ and $Z = \sum^r_{k= 3} Z_k$ under normal form and a smooth Hamiltonian $R$ satisfying $\partial^\alpha R(0)=0$ for all $\alpha\in \N^\Zc$ with $|\alpha|\geq r$, such that 
\begin{equation}\label{Y0}
(H_0 +  P) \circ \Phi_{\chi}^1 = H_0 + Z+R. 
\end{equation}
Then   the exponential estimate will by obtained by optimizing   the choice of $r$ and $N$.\\
We recall that for $\chi$ and $K$ to Hamiltonian, we have for all $k \geq 0$ 
$$
\frac{\dd^k }{\dd t^k} (K \circ \Phi_\chi^t )=    \{\chi, \{ \cdots \{\chi,K\}\cdot \}(\Phi_\chi^t)   =  (\mathrm{ad}^k_{\chi} K)(\Phi_\chi^t),
$$
where $\mathrm{ad}_{\chi}K = \{ \chi,K\}$. On the other hand, if $K$, $L$ are homogeneous polynomials of degree respectively $k$ and $\ell$ then $\{K,L\}$ is a homogeneous polynomial of degree $k+l-2$.  
Therefore, we obtain by using the Taylor formula
\begin{equation}\label{Y1}
(H_0+P) \circ \Phi_{\chi}^1  - (H_0+P) = \sum_{k = 0}^{r-3}
\frac{1}{(k+1)!} \mathrm{ad}^{k}_{\chi} (\{\chi,H_0+P\}) +\mathcal O_r 
\end{equation}
where $\mathcal O_r$ stands for any smooth function $R$ satisfying $\partial^\alpha R(0)=0$ for all $\alpha\in \N^\Zc$ with $|\alpha|\geq r$.
  Now we know that  for $\zeta\in\C$, the following relation holds:   
$$
\left(\sum_{k = 0}^{r-3} \frac{B_k}{k!} \zeta^k \right)\left(\sum_{k = 0}^{r-3} \frac{1}{(k+1)!} \zeta^k\right)=1+O(|\zeta|^{r-2})
$$
where $B_k$ are the Bernoulli numbers defined by the expansion of the generating function $\frac{z}{e^z - 1}$. Therefore, defining the two differential operators 
$$
A_r= \sum_{k = 0}^{r-3} \frac{1}{(k+1)!}\mathrm{ad}^{k}_{\chi} \quad \text{ and }\quad B_r=\sum_{k = 0}^{r-3} \frac{B_k}{k!} \mathrm{ad}^{k}_{\chi},
$$
we get
$$B_rA_r=\mathrm{Id}+C_r$$
where $C_r$ is a differential operator satisfying
$$
C_r \mathcal O_3= \mathcal O_r.$$
Applying $B_r$ to the two sides of equation \eqref{Y1}, we obtain
$$
\{\chi,H_0+P\}=B_r(Z-P)+\mathcal O_r.$$
Plugging the decompositions in homogeneous polynomials of $\chi$, $Z$ and $P$ in the last equation and equating the terms of same degree, we obtain after a straightforward calculus, the following recursive equations
\begin{equation}\label{Y2}
\{\chi_m,H_0\}-Z_m=Q_m,\quad m=3,\cdots,r,
\end{equation}
where
\begin{align}\begin{split}\label{Q_m}
Q_m&=-P_m + \sum_{k = 3}^{m-1} \{P_{m+2 - k},\chi_{k}\}\\
&+ \sum_{k =  1}^{m-3} \frac{B_k}{k!} \sum_{\substack{\ell_1+ \cdots \ell_{k+1} = m + 2k \\ 3\leq \ell_i \leq m - k} }
\mathrm{ad}_{\chi_{\ell_1}} \cdots \mathrm{ad}_{\chi_{\ell_{k}}}( Z_{\ell_{k+1}} - P_{\ell_{k+1}}).
\end{split}\end{align}
Notice that in the last sum, $\ell_i \leq m - k$ as a consequence of $3\leq \ell_i $ and $\ell_1+ \cdots \ell_{k+1} = m + 2k$.\\
Once these recursive equations solved, we define the remainder term as 
$R=(H_0+P) \circ \Phi_{\chi}^1 - H_0 - Z$. By construction, $R$ is analytic on a neighborhood of the origin in $\mathcal L_\rho$ and $R=\mathcal O_r$. As a consequence, by the Taylor formula,
 \begin{align}\begin{split}\label{R}
R&= \sum_{m\geq r+1}\sum_{k = 1}^{m-3} \frac{1}{k!} \sum_{\substack{\ell_1+ \cdots \ell_{k} = m + 2k \\ 3\leq \ell_i \leq r} }
\mathrm{ad}_{\chi_{\ell_1}} \cdots \mathrm{ad}_{\chi_{\ell_{k}}} H_0\\
&+\sum_{m\geq r+1} \sum_{k =  0}^{m-3} \frac{1}{k!} \sum_{\substack{\ell_1+ \cdots \ell_{k+1} = m + 2k \\ 3\leq \ell_1+ \cdots \ell_{k} \leq r\\ 3\leq\ell_{k+1}} }
\mathrm{ad}_{\chi_{\ell_1}} \cdots \mathrm{ad}_{\chi_{\ell_{k}}}P_{\ell_{k+1}}.
\end{split}\end{align}

\begin{lemma}\label{chim}
Assume that  the non resonance condition \eqref{Enonres1} is fulfilled. 
Let $r$ and $N$ be fixed. For $m=3,\cdots,r$, there exist homogeneous polynomials $\chi_m$ and $Z_m$ of degree $m$, with $Z_m$ in $N-$normal form, solutions of the recursive equation \eqref{Y2} and satisfying
\begin{equation}
\label{EchiZ}
\Norm{\chi_m}{} + \Norm{Z_m}{}  \leq (C m N^{\nu})^{m^2}
\end{equation}
where the constant $C$ does not depend on $r$ or $N$.
\end{lemma}
\begin{Proof}
We define $\chi_m$ and $Z_m$ by induction using Proposition \ref{Phomolo}. Note that \eqref{EchiZ} is clearly satisfied for $m = 3$, provided $C$ is big enough. Estimate \eqref{Ehomol}, together with \eqref{Ebrack} and the estimate on the Bernoulli numbers, $|B_k|\leq k!\ c^k$ for some $c>0$,  yields for all $m \geq 3$, 
 \begin{align*}\begin{split}
\gamma c_0^m N^{-\nu m}& \Norm{\chi_m}{} + \Norm{Z_m}{} \leq \Norm{P_m}{} + 2 \sum_{k = 3}^{m-1} k(m+2-k) \Norm{P_{m+2 - k}}{}\Norm{\chi_{k}}{} \\
&+ 2 \sum_{k =  1}^{m-3} (Cm)^k  \sum_{\substack{\ell_1+ \cdots \ell_{k+1} = m + 2k \\ 3\leq \ell_i \leq m - k} }
\ell_1\Norm{\chi_{\ell_1}}{}\cdots \ell_k\Norm{\chi_{\ell_{k}}}{} \Norm{ Z_{\ell_{k+1}} - P_{\ell_{k+1}}}{}.
\end{split}\end{align*} 
for some constant $C$. 
We set   $\beta_m = m (\Norm{\chi_m}{} + \Norm{Z_m}{})$.    Using $\Norm{P_m}{} \leq M R_0^{-m}$ (see end of subsection \ref{sec:NR}), we obtain  
 \begin{align*}\begin{split}
\beta_m &\leq  \beta^{(1)}_m+ \beta^{(2)}_m\quad \text{where}\\
\beta^{(1)}_m&=(C N^{\nu })^m m^3  \sum_{k = 3}^{m-1}   \beta_k  \quad \text{and}\\
\beta^{(2)}_m&=N^{\nu m }(Cm)^{m-1} \sum_{k =  1}^{m-3}   \sum_{\substack{\ell_1+ \cdots \ell_{k+1} = m + 2k \\ 3\leq\ell_i \leq m - k}} 
\beta_{\ell_1}{}\cdots \beta_{\ell_{k}}{} (\beta_{\ell_{k+1}} + \Norm{P_{\ell_{k+1}}}{})
\end{split}\end{align*}
where $C$  depends on $M$, $R_0$, $\gamma$ and $c_0$. 
It remains to prove by recurrence  that $\beta_m \leq (C m N^\nu)^{m^2}$, $m\geq 3$. Again this is true for $m=3$ adapting $C$ if necessary. Thus assume   that $\beta_j \leq (C j N^\nu)^{j^2}$  $j = 3,\ldots,m-1$, we then get for 
$$
\beta^{(1)}_m \leq (C N^{\nu })^m m^4  (C m N^\nu)^{(m-1)^2} 
\leq  (C m N^\nu)^{m^2 - m +1}
\leq \frac12 (C m N^\nu)^{m^2}
$$
as soon as $m \geq 4$, and provided $C > 2$.  
On the other hand, since $\Norm{P_m}{} \leq M R_0^{-m}$, we can  assume that   $\Norm{P_{\ell_{k+1}}}{} \leq \beta_{\ell_{k+1}}$    and we get 
$$
\beta^{(2)}_m
\leq N^{\nu m }(Cm)^{m-1} \sum_{k =  1}^{m-3}   \sum_{\substack{\ell_1+ \cdots \ell_{k+1} = m + 2k \\ 3\leq\ell_i \leq m - k}} 
(CN^\nu(m-k))^{\ell_1^2+\cdots +\ell_{k+1}^2} .
$$
Notice that the maximum of $\ell_1^2+\cdots +\ell_{k+1}^2$ when $\ell_1+ \cdots \ell_{k+1} = m + 2k $ and $ 3\leq\ell_i \leq m - k$ is obtained for $\ell_1=\cdots=\ell_k=3$ and $\ell_{k+1}= m-k$ and its value is $(m-k)^2+9k$. Furthermore
the cardinal of $\{\ell_1+ \cdots \ell_{k+1} = m + 2k, \ 3\leq\ell_i \leq m - k\}$ is smaller than $m^{k+1}$, hence
we obtain
$$
\beta^{(2)}_m
\leq \max_{k=1,\cdots,m-3}N^{\nu m }(Cm)^{m-1} Cm^{k+2}  (CN^\nu(m-k))^{(m-k)^2+9k} \leq \frac12 (C m N^\nu)^{m^2}
$$
for all $m\geq 4$ and adapting again $C$ if necessary.
\end{Proof}


\subsection{Normal form result}

For a number   $R_0$   , we set 
$  
B_\rho(R_0) = \{ z \in \Lc_\rho\, | \, \Norm{z}{\rho} < R_0  \}. 
$  
\begin{theorem}
\label{P51}
Assume that $P$ is analytic on a ball   $B_\rho(R_0)$ for some $R_0>0$   and $\rho>0$. Assume that the non resonance condition \eqref{Enonres1} is satisfied, and let $\beta < 1$ and $M>1$ be fixed. Then there exist constants $\varepsilon_0>0$ and $\sigma>0$  such that for all $\varepsilon < \varepsilon_0$,   there exists:  a polynomial $\chi$, a polynomial $Z$ in $|\ln\varepsilon|^{1 + \beta}$ normal form, and a   Hamiltonian $R$ analytic on $B_\rho(M\eps)$,     such that 
\begin{equation}
\label{FoNo}
(H_0 + P) \circ \Phi_\chi^1 = H_0 + Z + R.
\end{equation}
Furthermore, for all $z \in B_\rho(M\varepsilon)$, 
\begin{equation}
\label{Eesteps}
 \Norm{X_Z(z)}{\rho} \ + \Norm{X_\chi (z)}{\rho}  \leq 2\varepsilon^{3/2}, \quad \mbox{and}\quad \Norm{X_R (z)}{\rho}  \leq \eps\ e^{-\frac1 4 |\ln\varepsilon|^{1+\beta}}. 
\end{equation}
\end{theorem}

\begin{Proof}
Using Lemma \ref{chim}, for all $N$ and $r$, we can construct  polynomial Hamiltonians
$$
\chi(z) = \sum_{k = 3}^r \chi_k(z) \quad\mbox{and}\quad  Z(z) = \sum_{k = 3}^r Z_k(z) ,
$$
with $Z$ in $N$-normal form, such that \eqref{FoNo} holds with $R=\mathcal O_r$.
Now 
for fixed $\varepsilon>0$, we choose
$$
N\equiv N(\eps) = |\ln \varepsilon|^{1+\beta}\quad \mbox{and}\quad r\equiv r(\eps) = |\ln \varepsilon|^\beta. 
$$
This choice is motivated by the necessity of a balance between $Z$ and $R$ in \eqref{FoNo}: The error induced by $Z$ is controlled as in Remark \ref{choixN}, while the error induced by $R$ is controlled by Lemma \ref{chim}. 
By \eqref{EchiZ}, we have 
\begin{equation}
\label{W2}
\begin{array}{rcl}
\Norm{\chi_k}{} &\leq& (CkN^\nu)^{k^2} \leq \exp(k( \nu k(1+\beta)\ln |\ln \eps| + k \ln Ck))\\
&\leq& \exp(k( \nu r(1+\beta)\ln |\ln \eps| + r \ln Cr))\\
&\leq& \exp(k|\ln\eps|( \nu |\ln\eps|^{  \beta - 1}(1+\beta)\ln |\ln \eps| + |\ln\eps|^{  \beta - 1   } \ln C|\ln\eps|^{\beta}))\\
&\leq&\eps^{-k/8}, 
\end{array}
\end{equation}
  as $ \beta < 1$, and  for $\varepsilon \leq  \varepsilon_0$ sufficiently small.    
Therefore using Proposition \ref{P1}, we obtain for $z \in B_\rho(M\varepsilon)$
$$
|\chi_k(z)| \leq\eps^{-k/8}(M\eps)^k \leq  M^k\varepsilon^{7k/8}
$$
and thus
$$
|\chi(z)|\leq \sum_{k\geq 3}  M^k   \varepsilon^{7k/8}\leq \eps^{3/2}
$$
for $\eps$ small enough.
Similarly, we have for all $k \leq r$, 
$$ 
\Norm{ X_{\chi_k}(z)}{\rho}    \leq   2k\eps^{-k/8}(M\eps)^{k-1} \leq  2k M^{k-1}\varepsilon^{7k/8-1}   
$$
and
$$  
\Norm{X_\chi(z)}{\rho} \leq  \sum_{k\geq 3} 2 kM^{k-1}\varepsilon^{7k/8-1}\leq C\eps^{-1}\eps^{\frac{21}{8}}\leq \eps^{3/2}   
$$
for $\eps$ small enough.
Similar   bounds    clearly hold for $Z = \sum_{k = 3}^r Z_k$, which shows the first estimate in \eqref{Eesteps}. 

On the other hand, using $\mathrm{ad}_{\chi_{\ell_{k}}} H_0=Z_{\ell_{k}}+Q_{\ell_{k}}$   (see \eqref{Y2}), then using Lemma \ref{chim} and the definition of $Q_m$ (see \eqref{Q_m}),    we get $ \Norm{\mathrm{ad}_{\chi_{\ell_{k}}} H_0}{}\leq (CkN^\nu)^{{\ell_{k}}^2} \leq \eps^{-{\ell_{k}}/8}$, where the last inequality proceeds as in \eqref{W2}. Thus, using \eqref{R}, \eqref{W2}   and $\Norm{P_{\ell_{k+1}}}{}\leq MR_0^{-\ell_{k+1}}$   we obtain by Proposition \ref{P1} that for $z \in B_\rho(M\varepsilon)$
$$  
 \Norm{X_R(z)}{\rho}  \leq \sum_{m\geq r+1}\sum_{k=0}^{m-3}m(Cr)^{3m}\eps^{-\frac{m+2k}{8}}  \eps^{m-1}   \leq \sum_{m\geq r+1}m^2(Cr)^{3m}\eps^{m/2}   \leq (Cr)^{3r}\eps^{r/2}.
$$ 
Therefore, since $r=|\ln \eps|^\beta$, we get  $\Norm{X_R(z)}{\rho} \ \leq \eps\ e^{-\frac1 4 |\ln\varepsilon|^{1+\beta}}$ for $z \in B_\rho(M\varepsilon)$ and $\eps$ small enough.
\end{Proof}

\subsection{Bootstrap argument}
We are now in position to
 prove the main theorem of Section 1 which is actually a consequence of Theorem \ref{P51}. \\
Let $u_0\in\Ac_{2\rho}$ with $|u_0|_{2\rho}=\eps$ and denotes by $z(0)$ the corresponding sequence of its Fourier coefficients which belongs, by Lemma \ref{AL}, to in $\Lc_{\frac 3 2 \rho}$ with $||z(0)||_{\frac 3 2 \rho}\leq \frac{c_\rho}{4} \eps$ with   $c_\rho= \frac {2^{d+2} }{(1 - e^{-{\rho}/{2\sqrt{d}}})^d}$.   Let $z(t)$ be the local solution in $\Lc_\rho$  of the Hamiltonian system associated with $H=H_0+P$.\\
Let $\chi$, $Z$ and $R$ given by Theorem \ref{P51} with $M=c_\rho$ and let $y(t)=\Phi_\chi^1(z(t))$. We recall that since $\chi(z)=O(\Norm{z}{}^3)$, the transformation $\Phi_\chi^1$ is close to the identity, $\Phi_\chi^1(z)=z+O(\Norm{z}{}^2)$ and thus, for $\eps$ small enough,    we have  $\Norm{y(0)}{\frac 3 2 \rho}\leq \frac{c_\rho}{2} \eps$.
In particular, as noticed in \eqref{RN}, $R_N^\rho(y(0))\leq  \frac{c_\rho}{2} \eps\ e^{-\frac{\rho} 2 N}\leq \frac{c_\rho}{2} \eps\ e^{-\sigma N}$    where  $\sigma = \sigma_\rho \leq \frac{\rho}{2}$.\\   
Let $T_\eps$ be maximum of time $T$ such that $R_N^\rho(y(t))\leq  {c_\rho} \eps\ e^{-\sigma N}$ and $\Norm{y(t)}{\rho}\leq {c_\rho} \eps$ for all $|t|\leq T$. 
By construction, 
$$y(t)=y(0)+\int_0^t X_{H_0+Z}(y(s))ds +\int_0^t X_{R}(y(s))ds$$
 so using \eqref{ERN} for the first flow and \eqref{Eesteps} for the second one, we get for $|t|<T_\eps$,
\begin{align}\begin{split}\label{Z1}
R_N^\rho(y(t))&\leq  \frac 1 2 {c_\rho} \eps\ e^{-\sigma N}+   4   |t|\sum_{k=3}^r\Norm{Z_k}{} k^3(c_\rho \eps)^{k-1}e^{-2\sigma  N} +|t|\eps\ e^{-\frac1 4 |\ln\varepsilon|^{1+\beta}}\\
&\leq \left( \frac 1 2 +   4   |t|\sum_{k=3}^r\Norm{Z_k}{} k^3(c_\rho \eps)^{k-2}e^{-\sigma N}+|t|\eps\ e^{-\frac1 8 |\ln\varepsilon|^{1+\beta}}\right){c_\rho} \eps\ e^{-\sigma N}
\end{split}\end{align}
where in the last inequality we used $\sigma=\min \{ \frac 1 8, \frac{\rho} 2\}$ and $N=\ln |\eps|^{1+\beta}$.\\
Using Lemma \ref{chim}, we then verify
$$
R_N^\rho(y(t))\leq \left( \frac 1 2 +C|t|\eps\ e^{-\sigma N}\right){c_\rho} \eps\ e^{-\sigma N}
$$
and thus, for $\eps$ small enough,
\begin{equation}
\label{ca vient}
R_N^\rho(y(t))\leq {c_\rho} \eps\ e^{-\sigma N} \quad \text{for all}\quad |t|\leq \min \{T_\eps,e^{\sigma N}\}.
\end{equation}
Similarly we obtain
\begin{equation}
\label{Oh oui}
\Norm{y(t)}{\rho}\leq {c_\rho} \eps \quad \text{for all}\quad |t|\leq \min \{T_\eps,e^{\sigma N}\}.
\end{equation}
In view of the definition of $T_\eps$, \eqref{ca vient} and \eqref{Oh oui} imply $T_\eps\geq e^{\sigma N}$. In particular $\Norm{z(t)}{\rho}\leq 2{c_\rho} \eps$ for $|t|\leq e^{\sigma N}=\eps^{-\sigma |\ln \eps|^\beta}$ and using \eqref{Ediff2}, we finally obtain \eqref{mainest1} with    $C=\frac{2^{2d+5}}{(1-e^{-\rho/{2\sqrt{d}}})^{2d}}$.\\   
Estimate \eqref{mainest2} is an other consequence of the normal form result and Proposition \ref{Pcrux}. Actually we  use that the Fourier coefficients of $u(t)$ are given by $z(t)$ which is $\eps^2$-close to $y(t)$ which in turns is almost invariant:  in view of \eqref{marre} and as in \eqref{Z1}, we have 
$$
\sum_{j\in\Z}e^{\rho|j|}\big||y_j(t)|-|y_j(0)|\big|\leq \left(   4  |t|\sum_{k=3}^r\Norm{Z_k}{} k^3(c_\rho \eps)^{k-1}e^{-2\sigma N}+|t|\eps\ e^{-\frac1 4 |\ln\varepsilon|^{1+\beta}}\right)$$
from which we deduce
$$
\sum_{j\in\Z}e^{\rho|j|}\big||y_j(t)|-|y_j(0)|\big|\leq |t|\ e^{-\sigma N}
$$
and then \eqref{mainest2}.

\appendix
\section{Proof of the non resonance hypothesis}
Instead of proving Proposition \ref{res.1}, we prove a slightly more general result.
For a multi-index   $\jb\in \Zc^r$    we define
$$N(\jb)=\prod_{k=1}^{r}(1+|j_k|).$$

\begin{proposition} \label{res.2}
Fix $\gamma >0$ small enough and $m>d/2$.  There exist  positive
constants  $C$ and $\nu$ depending only on $m$, $R$ and $d$,  and a set $F_{\gamma} \subset \W_{m}$
whose measure is larger than $1- 4 \gamma$ such that if $V\in F_{\gamma}$
then for any $r\geq 1$
\begin{eqnarray}
\label{nr.d}
\left|\Omega(\jb)+\eps_{1} \omega_{\ell_{1}}+\eps_{2}
\omega_{\ell_{2}}\right|
\geq \frac{C^r \gamma^7}{N(\jb)^\alpha}
\end{eqnarray}
for any  $\jb \in\Zc^{r}$, 
for any indexes 
$l_{1},\ l_{2}\in \Z^d$,
and for any $\eps_{1},\eps_{2} \in \{ 0,1,-1\}$ such that $(\jb,(\ell_1,\eps_1),(\ell_2,\eps_2)) \notin \Nc_r$ is non resonant.

\end{proposition}

In order to prove proposition \ref{res.1}, we first prove that 
$\Omega(\jb)$ cannot accumulate on $\Z$. Precisely we have

\begin{lemma}
\label{lem.res.d}
Fix $\gamma >0$  and $m> d/2$.  There exist $0<C<1$ depending only on $m$, $R$ and $d$ and a set $F'_{\gamma} \subset \W_{m}$
whose measure is larger than $1-  4\gamma$ such that if $V\in F'_{\gamma}$
then for any $r\geq 1$
\begin{eqnarray}
\label{nr.d2}
\left|\Omega(\jb)-b\right|
\geq \frac{C^r\gamma}{N(\jb)^{m+d+3}}
\end{eqnarray}
for any non resonant $\jb \in\Z^{r}$ and for any $b\in \Z$.

\end{lemma}

\noindent
\begin{Proof}
Let  $(\alpha_{1},\ldots,\alpha_{r}) \neq 0$ in $\Z^r$, 
$M>0$ and 
$c\in \R$. By induction we can prove that the Lesbegue measure of 
$$
\{x\in[-M,M]^r \mid \va{\sum_{i=1}^r \alpha_{i}x_{i} +c}<\eta \}
$$
is smaller than $(2M)^{r-1} 2\eta $.  Hence given $\jb=(a_i,\delta_i)_{i=1}^r\in \Zc^r$,   and $b\in \Z$, the Lesbegue measure of
$$
\mathcal X_\eta:=\left\{ x\in[-1/2,1/2]^{r}\ :\ \left|\sum_{i=1}^r
\delta_i(|a_i|^2+x_i)-b \right|<\eta  \right\}
$$
is smaller than $2\eta$. Now consider the set
\begin{equation}
\label{saloperie}
\left\{V\in \W_m \mid \left|\Omega(\jb)-b\right|<\eta \right\}=\left\{V\in \W_m \mid  \left|\sum_{i=1}^r
\delta_i(|a_i|^2+\frac{v_{a_i}R}{(1+|a_i|)^m})-b \right|<\eta \right\}, 
\end{equation}
it is contained in the set of the $V$'s such that
$(Rv_{a_i}/(1+|a_i|)^m)_{i=1}^r\in\mathcal X_\eta$.  Hence the  measure of \eqref{saloperie} is smaller than $2R^{-r}N(\jb)^m \eta$. To conclude the proof we have
to sum over all the $\jb$'s and all the $b$'s. Now for a given $\jb$, remark that if
$|\Omega(\jb) -b|\geq \eta$ with $\eta \leq1$ then
$\va{b}\leq 2 N(\jb)^2$. So that  
to guarantee 
\eqref{nr.d2} for all possible choices of $\jb$,  $b$ and $r$, it 
suffices 
to remove from $\W_{m}$ a set of measure
$$
4\gamma\sum_{\jb\in \Zc^r} \frac{ C^r}{ R^rN(\jb)^{m+3+d}}N(\jb)^{m+2}  \leq 4\gamma \left[ \frac{2C}{R}\sum_{\ell \in \Z^d}\frac{1}{(1+|\ell|)^{d+1}}\right]^r .
$$
 Choosing $C \leq \frac12R\big(\sum_{\ell \in \Z^d}\frac{1}{(1 + |\ell|)^{d+1}}\big)^{-1}$ proves the result. \end{Proof}

\medskip

\noindent
\begin{Proofof}{proposition \ref{res.2}}
First of all, for $\eps_{1}=\eps_{2}=0$, \eqref{nr.d} is a direct consequence of lemma 
\ref{lem.res.d} choosing $\nu \geq m+d+3 $, $\gamma\leq 1$ and 
$F_{\gamma}=F'_{\gamma}$.\\
When $\eps_{1} =\pm 1$ and $\eps_{2}=0$, \eqref{nr.d} reads
\begin{equation}\label{l=1}
\va{\Omega(\jb) \pm \omega_{\ell_1}}\geq  \frac{C^r\gamma}{N(\jb)^\nu}.
\end{equation}
 Notice that $\va{\Omega(\jb) }\leq N(\jb)^2$ 
and thus, if $|\ell_1|\geq 2N(\jb)$, \eqref{l=1} is always true. 
When $|\ell_1|\leq 2 N(\jb)$, using that $N(\jb,\ell)=N(\jb)(1+|\ell_1|)$, we get applying lemma \ref{lem.res.d} with $b = 0$,
$$
\va{\Omega(\jb) + \varepsilon_1\omega_{\ell_1}}
= \va{\Omega(\jb,(\ell_1,\varepsilon_1))}\geq  \frac{C^{r+1}\gamma}{N(\jb)^{m+d+3}(3 N(\jb))^{m+d+3}}\geq 
\frac{\tilde C^{r}\gamma}{N(\jb)^\nu}
$$
with $\nu=2(m+d+3)$ and $\tilde C=\frac{C^2}{3^{m+d+3}}$. 
 In the 
same way we prove \eqref{nr.d} when $\eps_{1}\eps_{2}=1$ with the same choice of $\nu$. So it 
remains to establish an estimate of the form
\begin{equation}\label{l=2}
\va{\Omega(\delta,\jb)+ \omega_{\ell_{1}}-\omega_{\ell_{2}}}
\geq  \frac{\tilde C^r\gamma^4}{N(\jb)^\nu}.
\end{equation}
Assuming $|\ell_1| \leq |\ell_2|$, we have 
$$
|\omega_{\ell_{1}}-\omega_{\ell_{2}} - \ell_{1}^{2} +\ell_{2}^{2} |\leq 
\left|\frac{R|v_{\ell_{1}}|}
{(1+\va{\ell_{1}})^m}
-\frac{R|v_{\ell_{2}}|}{(1+\va{\ell_{2}})^m}\right|\leq \frac R {(1+\va{\ell_{1}})^m}.
$$
Therefore if ${(1+\va{\ell_{1}})^m}\geq \frac{2R}{C^r\gamma} N(\jb)^{m+d+3}$, we obtain \eqref{l=2} directly from lemma \ref{lem.res.d} applied with $b=\ell_{1}^{2} -\ell_{2}^{2}$ and choosing
$\nu =m+d+3$, $\tilde C = C/2$ and  $F_{\gamma}=F'_{\gamma}$.\\
Finally assume $(1+\va{\ell_{1}})^m\leq \frac{2R}{C^r\gamma} N(\jb)^{m+d+3}$, taking into acount 
$\va{\Omega(\jb) }\leq N(\jb)^2$, \eqref{l=2} is satisfied 
when $ \ell_{2}^{2} -\ell_{1}^{2} \leq 2 N(\jb)^2$. So it remains to consider the case  
when  $$1+|\ell_{1}| \leq 1+| \ell_{2}| \leq \left[\left(\frac{2R}{C^r\gamma} N(\jb)^{m+d+3}\right)^{2/m}+2 N(\jb)^2\right]^{1/2}\leq \left(\frac{3R}{C^r\gamma}\right)^\frac1 m N(\jb)^\frac{m+d+3}{m}.$$ 
Again we use lemma \ref{lem.res.d}  to conclude
\begin{align*}
\va{\Omega(\jb)+ \omega_{\ell_{1}}-\omega_{\ell_{2}}}&\geq
 \frac{C^{r+2}\gamma}{[N(\jb)(1+|\ell_{1}|)( 1+| \ell_{2}|)] ^{m+d+3}}\\&\geq
  \frac{C^{r+2}\gamma\left(\frac{C^r\gamma}{3R}\right)^\frac{m+d+3}{ m}}{N(\jb)^{m+d+3} N(\jb)^{2\frac{(m+d+3)^2}{m}}}\geq  \frac{\tilde C^r\gamma^{4+3/m}}{N(\jb)^\nu}
\end{align*}
with $\nu=m+d+3+(m+d+3)^2/m$ and $\tilde C=\frac{C^{(4m+d+3)/m}}{3R}.$
\end{Proofof}


\end{document}